\documentclass[pdflatex,sn-mathphys-num]{sn-jnl}

\usepackage{caption}
\usepackage{graphicx}%
\usepackage{multirow}%
\usepackage{amsmath,amssymb,amsfonts}%
\usepackage{amsthm}%
\usepackage{mathrsfs}%
\usepackage[title]{appendix}%
\usepackage{xcolor}%
\usepackage{textcomp}%
\usepackage{manyfoot}%
\usepackage{booktabs}%
\usepackage{algorithm}%
\usepackage{algorithmicx}%
\usepackage{algpseudocode}%
\usepackage{listings}%

\theoremstyle{thmstyleone}%
\newtheorem{theorem}{Theorem}
\newtheorem{lemma}[theorem]{Lemma}
\newtheorem{corollary}{Corollary}[theorem]
\theoremstyle{thmstyletwo}%
\newtheorem{example}{Example}%

\theoremstyle{thmstylethree}%
\newtheorem{definition}{Definition}%

\begin{document}

\title[FPT in LR-$(\delta,\phi)$ quasi partial $b$-metric spaces]{Common Fixed Point Theorems on LR-$(\delta,\phi)$ Quasi Partial $b$-Metric Spaces}


\author[1]{\fnm{} \sur{Anuradha Gupta}}\email{dishna2@yahoo.in}
\equalcont{These authors contributed equally to this work.}
\author*[2]{\fnm{Rahul} \sur{Mansotra}}\email{mansotrarahul2@gmail.com}
\equalcont{These authors contributed equally to this work.}
\affil[1]{\orgdiv{Department of Mathematics}, \orgname{ Delhi College of Arts and Commerce, University of Delhi}, \orgaddress{\street{Netaji Nagar}, \city{New Delhi}, \postcode{110023}, \state{Delhi}, \country{India}}}
\affil*[2]{\orgdiv{Department of Mathematics}, \orgname{University of Delhi}, \orgaddress{\street{} \city{New Delhi}, \postcode{110007}, \state{Delhi}, \country{India}}}

\abstract{In this article, we introduce the  LR-$(\delta,\phi)$ quasi partial $b$-metric space. Also, the common fixed point theorem on complete LR-$(\delta,\phi)$ quasi partial $b$-metric space has been proved. A non-trivial example is  also given.}




\maketitle

\vspace{-0.7cm}\section{Introduction and Preliminaries }\label{sec1}
Fixed point is an element which is invariant  under a self map. Fixed point theory involves finding the conditions under which the self map has a fixed point. Fixed point theory finds direction after the famous result of Banach Contraction Principle \cite{bm1}, which says every contraction on a complete metric space has a unique fixed point. Moreover, this principle gives a method for finding the fixed point in a complete metric space. From then, fixed point theorists trying to generalize this principle either by expanding the contraction condition or generalizing the domain space. The main motivation behind the  generalization of contraction condition is  uniform continuity  of contraction mapping  in the Banach Contraction Principle. In this direction, fixed point researchers works with the   generalization of  contraction conditions using  the terms   $d(x,Tx),d(x,y),d(y,Ty), d(y,Tx)$ and $d(x,Ty)$, where $x,y \in X,$ $d \text{ is a metric on } X,$ and $ T$ be a self-map on $X$. Few among  them are Kannan \cite{bm16}, Riech \cite{bm21}, Chatterjee  \cite{bm6}, Ciric \cite{bm22}   and Hardy-Rodgers \cite{bm11} contraction.  Rhoades \cite{bm22}, in 1977, did comparison among various contractive conditions in metric spaces. On the other hand, various generalizations of metric spaces is done by  modifying its properties to enhance the applicability of fixed point theorems in broader area. Some  of them are partial metric \cite{bm18},  dislocated metric \cite{bm23}, rectangular metric \cite{bm4}, $b$-metric \cite{bm7}, dislocated quasi-metric \cite{bm24} and quasi-partial $b$-metric space \cite{bm8}.
\par A point $x$ is called a coincidence point(common fixed point) of the pair $(U,V)$ if $Ux=Vx(Ux=Vx=x)$, where $U,V$ are two self-maps on a non-empty set $X$. Despite the fixed point theory for single valued mapping has numerous applications, significant number of pure and applied mathematics can be easily converted into more than one mapping problems. Common fixed points are the  representatives of  stability or equilibrium within the system. There are many real world scenarios like game theory and network theory, having  more than one player or variable, whom action can be taken as a mapping. Common fixed points give a stability point within them. Question of common fixed points was first raised by Isbell \cite{bm13} in the following problem:  
``Do any two continuous self maps on $[0,1]$ have a common fixed point."
 The answer is negative as discussed by Boyce \cite{bm3} and Huneke \cite{bm12} with the help of examples. Common fixed point theorems sometimes require certain assumptions on the  space or  mappings. Following, Jungck \cite{bm15} proposed  a definition of compatible mapping and gave common fixed point theorems for these mappings, which gives  direction in this area of study.   Many fixed point theorists trying to extends the general common fixed point theorems either by modifying  contraction or changing underlying space. For  more literature, we refer \cite{bm5, bm14, bm17, bm19, bm20, bm25, bm26}. \\
 Now, we will give the following prerequisites:
   \vspace{-0.4cm} \begin{definition}\cite{bm24}
      A mapping $q_{d}: X \times X \to [0,\infty)$ defines a dislocated quasi-metric on a non-empty set $X$ if for all $x,y,z \in X$, $q_{d}$ satisfies the following properties:
     \begin{itemize}
         \item[(i)] If $q_{d}(x,y)=q_{d}(y,x)=0$, then $x=y;$
         \item[(ii)] $q_{d}(x,y)\leq q_{d}(x,z)+q_{d}(z,y)-q_{d}(z,z).$
     \end{itemize}
     Then, the pair $(X,q_{d})$ is called dislocated quasi-metric space.
 \end{definition}
\vspace{-1cm}\begin{definition}\cite{bm8}
    A mapping $q_{pb} :X \times X \to [0,\infty)$ defines a quasi-partial $b$-metric  if there exists $s\geq 1$ such that for all $x,y, z \in X $,   $q_{pb}$ satisfies the following properties:
    \begin{enumerate}
        \item If $q_{pb}(x,y)=q_{pb}(x,x)=q_{pb}(y,y)$, then $x=y; $
        \item $q_{pb}(x,x) \leq q_{pb}(x,y)$;
        \item $q_{pb}(x,x) \leq q_{pb}(y,x)$;
        \item $ q_{pb}(x,y)\leq s(q_{pb}(x,z)+q_{pb}(z,y))-q_{pb}(z,z). $
    \end{enumerate}
    Then, $(X,q_{pb},s)$ is called quasi partial $b$-metric space with coefficient $s \geq 1$.
\end{definition}
\vspace{-1cm}\begin{definition}\cite{bm9}
    A closed  ball centered at $x$ with radius $\epsilon>0$ in a quasi-partial $b$-metric space $(X,q_{pb},s)$ is defined as follows: \begin{equation*}
        \overline{B^{l}}_{q_{pb}}(x, \epsilon)= \{ y \in X : q_{pb}(x,y) \leq \epsilon , q_{pb}(y,x) \leq \epsilon\}.
    \end{equation*} 
\end{definition}\vspace{-0.3cm} \noindent Now,  we define  the left closed ball in a quasi-partial $b$-metric space $(X,q_{pb},s)$ as follows: 
\begin{definition}
    A left closed  ball centered at $x \in X$ with radius $\epsilon>0$  in a quasi-partial $b$-metric space $(X,q_{pb},s)$ is  given as follows: 
    $$\overline{B^{l}}_{q_{pb}}(x, \epsilon)= \{ y \in X : q_{pb}(x,y) \leq \epsilon\}.$$ 
\end{definition}
\vspace{-1cm}\begin{definition}\cite{bm9}
   A subset of a quasi-partial $b$-metric space $(X,q_{pb},s)$ is said to be   closed set in $X$ if its complement in $X$ is open in $X$.
\end{definition}
\vspace{-0.5cm}\noindent It can be easily shown that a left closed ball in  a quasi-partial $b$-metric space is closed.
\vspace{-0.4cm}\begin{definition}\cite{bm2}
    let $s\geq 1$ be a fixed real number. A map $\psi : [0, \infty)  \to [0, \infty)$ is said to be a $b$-comparison function if $\psi$ is non-decreasing and $ \sum_{i=0}^{\infty}s^{i} \psi^{i}(t)$ converges for each $t\geq 0.$
\end{definition}
\vspace{-0.5cm} \noindent By $\Psi_{s}$, we denote the collection of all $b$-comparison functions.  Moreover, for each $\psi \in \Psi_{s}$, we can easily  find that $\psi(t)<t$, for each $t>0$.
\vspace{-0.4cm}\begin{lemma} \cite{bm8}
    Let $(X,q_{pb},s)$ be a quasi-partial $b$-metric space with  coefficient $s\geq 1$. Then the following statements holds: \begin{enumerate}
        \item If $q_{pb}(x,y)=0$, then $x=y.$
        \item If $x \neq  y$, then $q_{pb}(x,y)>0.$
    \end{enumerate} 
\end{lemma}
\vspace{-0.9cm}\begin{definition}\cite{bm24}
    Let $A \subseteq X$, $\delta : X \times X \to [0, \infty)$ and $T$ be a self map on $X$. Then  $T$ defines locally $\delta$-dominated mapping on $A$ if $\delta(x,Tx)\geq 1$,  for each $x \in A$. Moreover, the mapping $\delta$ defines locally triangular on $A$ if $\delta(x,y)\geq 1$ and $\delta(y,z)\geq 1$ implies $\delta(x,z)\geq 1$, for each $x,y,z \in A$.   \end{definition}
\vspace{-0.95cm}\begin{definition}\cite{bm24}
    Let $A \subseteq X$ and $\delta, \phi : X \times X \to [0, \infty)$ be two mappings. Then the pair $(\delta, \phi )$ is said to be locally  triangular on $A$ if $\delta (x,y) \geq \phi(x,y)$ and $\delta (y,z) \geq \phi(y,z)$ implies  $\delta (x,z) \geq \phi(x,z)$, for each $x,y,z \in A$. Moreover, a mapping $T: X \to X$ is known to be  locally $(\delta, \phi)$-dominated mapping on $A$ if $\delta(x,Tx)\geq \phi(x,Tx)$, for each $x \in A.$
\end{definition}
\vspace{-0.4cm}Shoaib et al. \cite{bm24} gave a concept of LR-complete dislocated quasi-metric space and showed that every  right-complete, left-complete and bi-complete dislocated quasi-metric space is LR-complete
dislocated quasi-metric space, but the converse may not be true. Shoaib and Mehmood \cite{bm24} 
introduced  LR-$(\delta,\phi)$ dislocated quasi-metric space and gave the common fixed point theorems on LR-$(\delta,\phi)$ dislocated quasi-metric spaces.
\par Motivated by the works of Shoaib and Mehmood \cite{bm24} on  LR-$(\delta,\phi)$ dislocated quasi-metric space and  inspired by the fixed point theorems on partial $b$-metric space by Gupta and Mansotra \cite{bm10}, we give the concept of LR-$(\delta,\phi)$ complete quasi-partial $b$-metric space and provide common fixed point theorem in the setting of LR-$(\delta,\phi)$ complete quasi-partial $b$-metric spaces.
\section{Common Fixed Point Theorem on LR-$(\delta,\phi)$ Quasi Partial $b$-Metric Space}\label{sec2}
We  start this section with the following definition:
\vspace{-0.45cm}\begin{definition}
    Let $(X,q_{pb},s)$ be a quasi-partial $b$-metric space with  coefficient $s\geq1$ and $\delta, \phi: X \times X \to [0, \infty)$ be  two mappings. Then, the  sequence $(x_{n})$  in $X$ with $\delta(x_{n},x_{n+1})\geq \phi(x_{n},x_{n+1})$ is said to be
    \begin{itemize}
        \item[(i)] $(\delta,\phi)$-Cauchy in $X$ if  $\lim \limits_{n,m\to \infty}q_{pb}(x_{m},x_{n})$ and $\lim \limits_{n,m\to \infty}q_{pb}(x_{n},x_{m})$ exist and finite.
        \item[(ii)] $(\delta,\phi)$-left or right converges to $x$ in $ X$ if $\lim \limits_{n \to \infty}q_{pb}(x_{n},x)=q_{pb}(x,x)$  or $\lim \limits_{n \to \infty}q_{pb}(x,x_{n})=q_{pb}(x,x).$   
         \item[(iii)] $(\delta,\phi)$-converges to $x$ in $X$ if it is $(\delta,\phi)$-left and right converges to $x$ in $X$.
        \item[(iv)] $X$ is LR-$(\delta,\phi)$ complete if each    $(\delta,\phi)$-Cauchy in $X$ is $(\delta,\phi)$-left or right convergent in $X$ i.e.\\
          $\lim \limits_{n,m\to \infty}q_{pb}(x_{m},x_{n})=\lim \limits_{n,m\to \infty}q_{pb}(x_{n},x_{m})=\lim \limits_{n \to \infty}q_{pb}(x_{n},x)=q_{pb}(x,x)$\\
         or \\
        $\lim \limits_{n,m\to \infty}q_{pb}(x_{m},x_{n})=\lim \limits_{n,m\to \infty}q_{pb}(x_{n},x_{m})=\lim \limits_{n \to \infty}q_{pb}(x,x_{n})=q_{pb}(x,x).$
        \end{itemize}
\end{definition}
\vspace{-0.9cm}\begin{theorem}\label{thm1}
Let $x_{0} \in X, \epsilon > 0, $ and $\delta, \phi : X \times X \to [0,\infty)$ be two mappings such that the pair $(\delta, \phi)$ is locally triangular  on left closed ball $\overline{B^{l}}_{q_{pb}}(x_{o}, \epsilon)$. Also,  $(X,q_{pb}, s)$ be an LR-$(\delta,\phi)$ 
 complete quasi-partial $b$-metric space with coefficient $s\geq 1$ and $U,V: X \to X $ be two locally $(\delta, \phi)$-dominated  self mappings  on $\overline{B^{l}}_{q_{pb}}(x_{o}, \epsilon)$. Suppose that 
 there exists a function $\psi \in \Psi_{s}$  such that the pair $(U,V)$  having the following properties: $\hspace{1cm}$\\
(i) For any $x,y \in \overline{B^{l}}_{q_{pb}}(x_{o}, \epsilon)$ such that either $\delta(x,y) \geq \phi(x,y) $ or $\delta(y,x) \geq \phi(y,x) $
     \begin{equation}
         \text{implies } \max \{ q_{pb}(Ux,Vy),  q_{pb}(Vy,Ux) \} \leq \psi (M_{s}(x,y)),
     \end{equation}
     where
     \begin{equation*}
        \hspace{-0.5cm} M_{s}(x,y)= \max
\Biggr \{ q_{pb}(x,y),    q_{pb}(x,Ux), q_{pb}(y,Vy), \frac{q_{pb}(x,Vy)+q_{pb}(y,Ux)-q_{pb}(x,x)}{2s}\Biggr \}.\end{equation*}
  \begin{align}
     \text{(ii) } q_{pb}(Vy,Ux) \leq q_{pb}(x,y),  \text{ for all } x,y \in \overline{B^{l}}_{q_{pb}}(x_{o}, \epsilon).\hspace{4.55cm}\end{align}
 \begin{equation}
   \text{(iii) } \sum_{i=0}^{j}s^{i+1} \psi^{i} ( \max  \{ q_{pb}(x_{0}, Ux_{0}), q_{pb}(Ux_{0},x_{0})\} ) \leq \epsilon, \text{ for all 
 } j \in \{ 0\}\cup \mathbb{N}.\hspace{1.8cm}
    \end{equation}
  Then $U$ and $V$ have atleast one common fixed point. Moreover, the common fixed point  set is singleton if any pair $(x,y)$ from the set of  common fixed points of   pair $(U,V)$ satisfies $\delta(x,y)\geq \phi(x,y).$
\end{theorem}
\vspace{-0.6cm}\begin{proof}
    Define a sequence $(x_{n})$ as $x_{2n+1}=Ux_{2n}$ and $x_{2n+2}=Vx_{2n+1}$ for each $n \in \mathbb{N} \cup \{ 0\}$. Since $x_{0} \in \overline{B^{l}}_{q_{pb}}(x_{o}, \epsilon)$ and $U$ is a locally $(\delta, \phi)$-dominated mapping on $\overline{B^{l}}_{q_{pb}}(x_{o}, \epsilon)$, we have
$$\delta(x_{0},x_{1})= \delta (x_{0}, Ux_{0}) \geq \phi(x_{0}, Ux_{0})=\phi(x_{0}, x_{1}).$$
Also, consider \begin{align*}
    q_{pb}(x_{0},x_{1})&= q_{pb}(x_{0},Ux_{0})\\
    & \leq sq_{pb}(x_{0}, Ux_{0}) + s^{2}\psi_{s}^{1} (q_{pb}(x_{0}, Ux_{0}) ) +... + s^{j+1}\psi_{s}^{j}q_{pb}(x_{0}, Ux_{0})\\&= \sum_{i=0}^{j}s^{i+1} \psi^{i} ( q_{pb}(x_{0}, Ux_{0}) ) \leq \sum_{i=0}^{j}s^{i+1} \psi^{i} ( \max \{ q_{pb}(x_{0}, Ux_{0}),  q_{pb}(Ux_{0}, x_{0})\} ) \\
    & \leq  \epsilon  , \hspace{7cm}\text{ by }
(3)\end{align*}
which implies that $x_{1} \in \overline{B^{l}}_{q_{pb}}(x_{o}, \epsilon).$ As $U$ is a locally $(\delta, \phi)$-dominated mapping  on $\overline{B^{l}}_{q_{pb}}(x_{o}, \epsilon)$, we get 
$$\delta(x_{1},x_{2})=\delta(x_{1},Vx_{1})\geq \phi(x_{1},Vx_{1})=\phi(x_{1},x_{2}). $$ We will proceed by the mathematical induction. Assume that $x_{2}, x_{3},..., x_{i}\in \overline{B^{l}}_{q_{pb}}(x_{o}, \epsilon)$. If $i$ is even, then $i=2j+2$, where $j \in \{ 0,1,2,...,\dfrac{i-2}{2}\}.$ As, $U$ is a locally $(\delta,\phi)$-dominated mapping on $ \overline{B^{l}}_{q_{pb}}(x_{o}, \epsilon)$, we get
\begin{equation}\delta(x_{2j+2},x_{2j+3})=\delta(x_{2j+2},Ux_{2j+2})\geq \phi(x_{2j+2},Ux_{2j+2})=\phi(x_{2j+2},x_{2j+3}) \end{equation} where $j \in \{ 0,1,2,...,\dfrac{i-2}{2}\}.$ As $x_{2j+1} \in \overline{B^{l}}_{q_{pb}}(x_{o}, \epsilon)$ for all $j \in \{ 0,1,2,...,\dfrac{i-2}{2}\}$ and $V $ is a locally $(\delta,\phi)$-dominated mapping, we get 
\begin{equation}
    \delta(x_{2j+1},x_{2j+2}) = \delta(x_{2j+1},Vx_{2j+1})\geq \phi (x_{2j+1},Vx_{2j+1})= \phi (x_{2j+1},x_{2j+2}),
\end{equation} for all $j \in \{ 0,1,2,...,\dfrac{i-2}{2}\}$. Let $q_{2k}=q_{pb}(x_{2k},x_{2k+1}), q^{'}_{2k}=q_{pb}(x_{2k+1},x_{2k})$ and $q^{''}_{2k}=q_{pb}(x_{2k},x_{2k})$. If either $q_{2j^{'}+1}$ or $q^{'}_{2j^{'}+1}=0$, for some $j^{'} \in \{ 0,1,2,...,\dfrac{i-2}{2}\}.$   Then by lemma $1$, we get $x_{2j^{'}+1}=x_{2j^{'}+2}$, which implies that $x_{2j^{'}+1}=Vx_{2j^{'}+1}$ and $x_{2j^{'}+2}=Vx_{2j^{'}+2}$. If either $q_{2j^{'}+2}$ or $q^{'}_{2j^{'}+2} =0$, then by Lemma 1, $x_{2j^{'}+2}=x_{2j^{'}+3}$. Hence, $x_{2j^{'}+2}$ is the common fixed point of the pair $(U,V)$. On the other hand, if $q_{2j^{'}+2}$ and $q^{'}_{2j^{'}+2} $ are not equal to zero, then by inequality $(5)$ and $(1)$, we have
\begin{align*}
    \max \{ q_{2j^{'}+2}, q^{'}_{2j^{'}+2}\}&= \max \{ q_{pb}(Ux_{2j^{'}+2},Vx_{2j^{'}+1}),q_{pb}(Vx_{2j^{'}+1},Ux_{2j^{'}+2}) \}\\&\leq  \psi \Biggr( \max \Biggr \{ q^{'}_{2j^{'}+1}, q_{pb}(x_{2j^{'}+2}, Ux_{2j^{'}+2}), q_{pb}(x_{2j^{'}+1}, Vx_{2j^{'}+1}),\\ & \hspace{1cm}\frac{q_{pb}(x_{2j^{'}+2},Vx_{2j^{'}+1})+q_{pb}(x_{2j^{'}+1},Ux_{2j^{'}+2})-q^{''}_{2k+2}}{2s} \Biggr \} \Biggr)\\&\leq  \psi \Biggr( \max \Biggr \{ q^{'}_{2j^{'}+1}, q_{2j^{'}+2}, q_{2j^{'}+1}, \frac{s(q_{2j^{'}+1}+q_{2j^{'}+2})-q^{''}_{2k+2}}{2s}\Biggr\} \Biggr)\\& \leq \psi \Biggr(\max \Biggr \{0,q_{2j^{'}+2},0, \frac{s(0+q_{2j^{'}+2})-q^{''}_{2k+2}}{2s}  \Biggr \} \Biggr)= \psi (q_{2j^{'}+2})\\
    & \leq \psi ( \max \{  q^{'}_{2j^{'}+2},  q_{2j^{'}+2} \})
\end{align*}
which is contrary to $\psi(t)< t, $ for each $t>0.$ Hence, $q^{'}_{2j^{'}+2} = q_{2j^{'}+2}=0.$\\
Assume that  $q_{2j+1}$  and  $q^{'}_{2j+1}$ are not equal to zero, for all $j \in \{ 0,1,2,...,\dfrac{i-2}{2}\}.$ By inequality $(4)$ and $(1)$, we have \begin{equation} \begin{split}
    \max \{  q^{'}_{2j+1}, q_{2j+1} \} &= \max \{ q_{pb}(Ux_{2j},Vx_{2j+1}), q_{pb}(Vx_{2j+1},Ux_{2j})  \}\\
    & \leq \psi \Biggr( \max \Biggr \{ q_{2j} , q_{pb}(x_{2j},Ux_{2j}), q_{pb}(x_{2j+1},Vx_{2j+1}), \\& \hspace{1cm} \frac  {q_{pb}(x_{2j},Vx_{2j+1}) + q_{pb}(x_{2j+1},Ux_{2j}) - q^{''}_{2j}}{2s} \Biggr \}  \Biggr) \\
    & \leq \psi \Biggr( \max \Biggr \{ q_{2j}, q_{2j+1},\frac  {s(q_{2j}+q_{2j+1}) - q^{''}_{2j}}{2s} \Biggr \}  \Biggr)\\&= \psi ( \max \{ q_{2j}, q_{2j+1} \} ). \end{split} 
\end{equation}If  $\max \{ q_{2j}, q_{2j+1} \} = q_{2j+1},$ then  by inequality $(6)$, we have
$$\max \{  q^{'}_{2j+1}, q_{2j+1} \} \leq  \psi ( q_{2j+1} )\leq\psi ( \max \{ q^{'}_{2j}, q_{2j+1} \} ),$$  which is contrary to $\psi(t)<t$, for each $t>0$. Hence,   $\max \{ q_{2j}, q_{2j+1} \} = q_{2j}$. Thus, by inequality $(6)$, we get \begin{equation}
    \max \{  q^{'}_{2j+1}, q_{2j+1} \} \leq \psi (q_{2j}) \leq \psi ( \max \{ q_{2j}, q^{'}_{2j}  \})
\end{equation}
Also, $\delta(x_{2j-1}, x_{2j}) \geq \phi(x_{2j-1}, x_{2j})$  for each  $j \in \{ 0,1,2,...,\dfrac{i-2}{2}\}$. So, by inequality $(1)$, we have
\begin{equation}\begin{split}
     \max \{ q_{2j}, q^{'}_{2j}  \}&= \max \{  q_{pb}(Ux_{2j},Vx_{2j-1}),q_{pb}(Vx_{2j-1},Ux_{2j})\}\\& \leq \psi \Biggr( \max \Biggr\{ q^{'}_{2j-1},   q_{pb}(x_{2j},Ux_{2j}), q_{pb}(x_{2j-1},Vx_{2j-1}),\\& \hspace{1cm} \frac{q_{pb}(x_{2j},Vx_{2j-1})+q_{pb}(x_{2j-1},Ux_{2j})-q^{''}_{2j}}{2s} \Biggr\} \Biggr)\\& \leq \psi \Biggr( \max \Biggr\{ q^{'}_{2j-1},   q_{2j}, q_{2j-1}, \frac{s(q_{2j-1}+q_{2j})-q^{''}_{2j}}{2s} \Biggr\} \Biggr)\\&= \psi ( \max \{ q^{'}_{2j-1},   q_{2j}, q_{2j-1} \} ).
\end{split}   
\end{equation}
If $\max \{ q^{'}_{2j-1},   q_{2j}, q_{2j-1} \}=q_{2j}, $ then by inequality $(8)$, we have 
$$
  \max \{ q_{2j}, q^{'}_{2j}  \}   \leq \psi (   q_{2j}) \leq \psi( \max \{ q_{2j}, q^{'}_{2j} \}),
$$ which is a    contradiction to   $\psi(t)<t,$ for each $t>0.$ Hence, $\max \{ q^{'}_{2j-1},   q_{2j}, q_{2j-1} \}=\max \{ q^{'}_{2j-1}, q_{2j-1} \}.$ Thus, by inequality $(8)$, we have
\begin{equation}
     \max \{ q_{2j}, q^{'}_{2j}  \} \leq \psi( \max \{ q^{'}_{2j-1}, q_{2j-1} \}).
\end{equation}
Also, $\psi$ is increasing. Hence, \begin{equation}
     \psi( \max \{ q_{2j}, q^{'}_{2j}  \} )\leq \psi^{2}( \max \{ q^{'}_{2j-1}, q_{2j-1} \}).
\end{equation}
Now, on using inequality $(10)$ in inequality $(7)$, we get 
\begin{equation}
   \max \{  q^{'}_{2j+1}, q_{2j+1} \}   \leq \psi^{2}( \max \{ q^{'}_{2j-1}, q_{2j-1} \}).
\end{equation} 
Also,  $\delta(x_{2j-2},x_{2j-1})\geq \phi(x_{2j-2},x_{2j-1})$, for each $j \in \{ 0,1,2,...,\dfrac{i-2}{2}\}$. Thus, by inequality $(1)$, we have 
\begin{equation}\begin{split}
     \max \{ q^{'}_{2j-1}, q_{2j-1} \}&=\max \{  q_{pb}(Ux_{2j-2},Vx_{2j-1}),q_{pb}(Vx_{2j-1},Ux_{2j-2})\}\\& \leq \psi \Biggr( \max \Biggr\{ q_{2j-2},   q_{pb}(x_{2j-2},Ux_{2j-2}), q_{pb}(x_{2j-1},Vx_{2j-1}),\\& \hspace{1cm} \frac{q_{pb}(x_{2j-2},Vx_{2j-1})+q_{pb}(x_{2j-1},Ux_{2j-2})-q^{''}_{2j-2}}{2s} \Biggr\} \Biggr)\\& 
     \leq \psi \Biggr( \max \Biggr\{ q_{2j-2}, q_{2j-1},\frac{s(q_{2j-2} + q_{2j-1})-q^{''}_{2j-2}}{2s} \Biggr\} \Biggr)\\& = \psi(\max \{ q_{2j-2}, q_{2j-1}\})
\end{split}  
\end{equation}
If $\max \{ q_{2j-2}, q_{2j-1}\}=q_{2j-1}$, then by inequality $(12)$, we have
\begin{equation*}
    \max \{  q^{'}_{2j-1}, q_{2j-1} \} \leq \psi (  q_{2j-1}) \leq \psi(  \max \{  q^{'}_{2j-1}, q_{2j-1} \}  ).
\end{equation*} 
which is contradiction to $\psi(t)<0$, for each $t  >0$. Hence, $$\max \{ q_{2j-2}, q_{2j-1}\}=q_{2j-2}.$$
Then, by inequality $(12)$, we get  \begin{equation}
     \max \{  q^{'}_{2j-1}, q_{2j-1} \} \leq \psi (  q_{2j-2}) \leq \psi(  \max \{  q^{'}_{2j-2}, q_{2j-2} \}  ).
\end{equation}
On substituting the inequality $(13)$ in  inequality $(11)$, we get 
\begin{equation*}
       \max \{  q^{'}_{2j+1}, q_{2j+1} \}   \leq \psi^{3}( \max \{  q^{'}_{2j-2}, q_{2j-2} \}).
\end{equation*}
On repeating the process, we get
\begin{equation}
     \max \{  q^{'}_{2j+1}, q_{2j+1} \}  \leq \psi^{2j+1}(  \max \{  q^{'}_{0}, q_{0} \}  ).
\end{equation}
If $i$ is odd, then $i=2j+1$, where $i \in \{  1,2,3,..., \frac{i-1}{2}\}$. Further, by  following the same steps, we get \begin{equation}
     \max \{  q^{'}_{2j}, q_{2j} \} \leq \psi^{2j}(  \max \{  q^{'}_{0}, q_{0} \}  ).
\end{equation} 
Combining  the two inequalities $(14)$ and $(15)$, we get 
\begin{equation}
     \max \{  q^{'}_{n}, q_{n} \} \leq \psi^{n}(  \max \{  q^{'}_{0}, q_{0} \}  ), \text{ for all } n \in \{  1,2, ..., i \}.
\end{equation}
Now, we claim that $x_{i+1} \in  \overline{B^{l}}_{q_{pb}}(x_{o}, \epsilon) $. For this, consider \begin{equation*}
    \begin{split}
        q_{pb}(x_{0},x_{i+1}) &\leq  sq_{0}+ s^{2}q_{1} +...+s^{i}q_{i-1}+s^{i}q_{i}-q^{''}_{1} -q^{''}_{2}-...-q^{''}_{j-1}-q^{''}_{j}\\& \leq s\max \{ q_{0}, q^{'}_{0}\}+ s^{2}\max \{ q_{1}, q^{'}_{1}\} +...+s^{i}\max \{ q_{i-1}, q^{'}_{i-1}\}+s^{i+1}\max \{ q_{i}, q^{'}_{i}\}\\& \leq   \sum_{n=0}^{i}s^{n+1} \psi^{n} ( \max  \{ q^{'}_{0}, q_{0}\} )\hspace{4.4cm}\text{ by  }(16)\\&
        \leq \epsilon. \hspace{8cm}\text{ by }(3)
    \end{split}
\end{equation*}
Hence,  $x_{i+1} \in  \overline{B^{l}}_{q_{pb}}(x_{o}, \epsilon),$ which further implies that $x_{n} \in \overline{B^{l}}_{q_{pb}}(x_{o}, \epsilon)$, for all $n \in \{0\} \cup \mathbb{N}$. We will show that $(x_{n}) $ is a LR-$(\delta,\phi)$ Cauchy sequence. Let $\epsilon > 0$, then there exists $k \in \mathbb{N}$, such that \begin{equation}
    \sum_{j\geq k}s^{j+1} \psi^{j} ( \max  \{ q^{'}_{0}, q_{0})\} ) < \epsilon.
\end{equation}
For $m>n\geq k$, consider 
\begin{equation*} \begin{split}
    q_{pb}(x_{m},x_{n}) &\leq sq^{'}_{n}+ s^{2}q^{'}_{n+1} +...+s^{m-n-1}q^{'}_{m-2}+s^{m-n-1}q^{'}_{m-1}-q^{''}_{n} -q^{''}_{n+1}-...-q^{''}_{m-1}\\& \leq s^{n+1}(\max \{q^{'}_{n}, q_{n} \})+ s^{n+2}(\max \{ q^{'}_{n+1}, q_{n+1} \} )+...\\&\hspace{1.5cm}+s^{m-1}(\max \{q^{'}_{m-2}, q_{m-2} \})+s^{m}(\max \{ q^{'}_{m-1},q_{m-1} \})\\& \leq   \sum_{j=n}^{m-1}s^{j+1} \psi^{j} ( \max  \{ q^{'}_{0}, q_{0}\} )\hspace{4.4cm}\text{ by }(16)\\& \leq   \sum_{j\geq k}s^{j+1} \psi^{j} ( \max  \{ q^{'}_{0}, q_{0}\} )\\& < \epsilon. \hspace{8cm}\text{ by }(17)
    \end{split}
\end{equation*} Thus, $(x_{n})$ is a right $(\delta, \phi)$-Cauchy sequence in $ \overline{B^{l}}_{q_{pb}}(x_{o}, \epsilon).$ In the similar way, we can easily verify that $(x_{n})$ is a left $(\delta, \phi)$-Cauchy sequence in  $ \overline{B^{l}}_{q_{pb}}(x_{o}, \epsilon)$. Thus, $(x_{n})$ is a LR-$(\delta, \phi)$ Cauchy sequence in  $ \overline{B^{l}}_{q_{pb}}(x_{o}, \epsilon)$. As, $(X,q_{pb},s)$ is a LR-$(\delta,\phi)$ complete quasi- partial $b$-metric space and every left closed ball of a LR-$(\delta, \phi)$ complete quasi-partial $b$-metric space is closed. Thus, $ \overline{B^{l}}_{q_{pb}}(x_{o}, \epsilon)$ is a     LR-$(\delta, \phi)$ complete. So, there is  $x \in \overline{B^{l}}_{q_{pb}}(x_{o}, \epsilon)$ such that $(x_{n})$ converges to $x$ i.e. $\lim \limits_{m,n \to \infty}q_{pb}(x_{n},x_{m})=\lim \limits_{m,n \to \infty}q_{pb}(x_{m},x_{n})=\lim \limits_{n \to \infty}q_{pb}(x_{n},x)=q_{pb}(x,x)=0$ or $\lim \limits_{m,n \to \infty}q_{pb}(x_{n},x_{m})=\lim \limits_{m,n \to \infty}q_{pb}(x_{m},x_{n})=\lim \limits_{n \to \infty}q_{pb}(x,x_{n})=q_{pb}(x,x)=0.$ 
\begin{equation*}
   \hspace{-5.1cm}\text{Also, } q_{pb}(Vx,Ux) \leq q_{pb}(x,x)=0, \hspace{2cm}\text{ by } (2)
\end{equation*}
Hence, by Lemma $1$, $Ux=Vx$. If $\lim \limits_{m,n \to \infty}q_{pb}(x_{n},x_{m})=\lim \limits_{m,n \to \infty}q_{pb}(x_{m},x_{n})=\lim \limits_{n \to \infty}q_{pb}(x_{n},x)=q_{pb}(x,x)=0$ holds, then \begin{equation*}\begin{split}
 q_{pb}(Vx,x) &\leq s(q_{pb}(Vx,Ux_{2n})+q_{pb}(Ux_{2n},x))- q_{pb}(Ux_{2n}, Ux_{2n})\\ & \leq  s(q_{pb}(x_{2n}, x)+q_{pb}(x_{2n+1},x)) \hspace{3cm}\text{ by  } (2)\end{split}
\end{equation*}
On applying  the limit, we get $ q_{pb}(Vx,x)=0$. Hence, by Lemma $1$, $Vx=x.$ 
If $\lim \limits_{m,n \to \infty}q_{pb}(x_{n},x_{m})=\lim \limits_{m,n \to \infty}q_{pb}(x_{m},x_{n})=\lim \limits_{n \to \infty}q_{pb}(x,x_{n})=q_{pb}(x,x)=0$ holds, then \begin{equation*}\begin{split}
 q_{pb}(x,Ux) &\leq s(q_{pb}(x,Vx_{2n+1})+q_{pb}(Vx_{2n+1},Ux))-q_{pb}(Vx_{2n+1},Vx_{2n+1}) \\ &\leq s(q_{pb}(x,x_{2n+2})+q_{pb}(x,x_{2n+1}))\hspace{3cm}\text{ by  } (2)\end{split}
\end{equation*}
On applying  the limit, we get $ q_{pb}(x,Ux)=0$. Hence, by Lemma $1$, $Ux=x.$ 
Thus, $x$ becomes the common fixed point of pair $(U,V)$ in $ \overline{B^{l}}_{q_{pb}}(x_{o}, \epsilon) $ . Now, we will show the uniqueness  of $x$ if $\delta(x,y)\geq \phi(x,y)$ for any other common fixed point $y$ in $ \overline{B^{l}}_{q_{pb}}(x_{o}, \epsilon) $. Let $y$ be another common fixed point  of  pair $(U,V)$ in $ \overline{B^{l}}_{q_{pb}}(x_{o}, \epsilon) $ such that  $\delta(x,y) \geq  \phi(x, y)$.  Thus, by inequality $(1)$, we have \begin{equation*}
\begin{split}
  \max \{ q_{pb}(x,y),  q_{pb}(y,x) \} &=\max \{ q_{pb}(Ux,Vy),  q_{pb}(Vy,Ux) \} \\&\leq \psi ( \max
\Biggr \{ q_{pb}(x,y),    q_{pb}(x,Ux), q_{pb}(y,Vy), \\ & \hspace{1cm} \frac{q_{pb}(x,Vy)+q_{pb}(y,Ux)-q_{pb}(x,x)}{2s}\Biggr \})\\& = \psi ( \max
\Biggr \{ q_{pb}(x,y),    q_{pb}(x,x), q_{pb}(y,y), \end{split}\end{equation*}\begin{equation*}
\begin{split}  & \hspace{4cm} \frac{q_{pb}(x,y)+q_{pb}(y,x)-q_{pb}(x,x)}{2s}\Biggr \}) \\&  \hspace{3cm}\leq \psi ( \max
 \{ q_{pb}(x,y),    q_{pb}(y,x) \}),
\end{split}
\end{equation*}
which implies that $x=y.$ 
\end{proof}
\vspace{-0.5cm}\begin{example}
    Let $X= [0,5)$. Define $q_{pb}: X \times X \to [0,\infty)$ as: \begin{center}
    $
q_{pb}(x,y) = 
\begin{cases}
x,      &   \text{if }x,y \in A,x=y \\
3x,     & \text{if } x \in B, y \in A \\
|x-y|^{2}, & \text{if }  x,y \in B \\
  \max \{ 2x,y\} +x, & \text{otherwise } 
\end{cases}
$\end{center} where $A=[0,4]$ and $B=(4,5)$.
We will show that $(X,q_{pb},s)$ is a quasi-partial $b$-metric space with coefficient $s=2$. For that we will claim the following:
\begin{itemize}
    \item[(i)]If $q_{pb}(x,y)=q_{pb}(x,x)=q_{pb}(y,y)$, then $x=y.$
\item[(ii)] $q_{pb}(x,x) \leq q_{pb}(x,y)$ and $q_{pb}(x,x) \leq q_{pb}(y,x)$.
\item[(iii)] $q_{pb}(x,y)\leq 2(q_{pb}(x,z)+q_{pb}(z,y))-q_{pb}(z,z).$
 \end{itemize}
 \noindent This can be easily  verified by the Table 1, Table 2, Table 3 respectively.

\begin{minipage}{0.3\textwidth}
\centering
\begin{tabular}{|c|c|c|}
\hline
$q_{pb}(x,y)$ & $q_{pb}(x,x)$ & $q_{pb}(y,y)$ \\ \hline
\multicolumn{3}{|c|}{$\text{Case I: For } x,y\in A, x \neq y $} \\ \hline
$\max \{ 2x, y \}+x$ & $x$ & $y$ \\ \hline
\multicolumn{3}{|c|}{$\text{Case II: For }x,y \in B, x \neq y$} \\ \hline
$|x-y|^{2}$ & $0$ & $0$ \\ \hline
\multicolumn{3}{|c|}{$\text{Case III: For }x \in A,y \in B$} \\ \hline
$\max \{ 2x, y \}+x$ & $x$ & $0$ \\ \hline
\multicolumn{3}{|c|}{$\text{Case IV: For }x \in B,y \in A$} \\ \hline
$3x$ & $0$ & $y$ \\ \hline
\end{tabular}\vspace{0.3cm}
\captionof{table}{}
\end{minipage}
\hfill
\begin{minipage}{0.48\textwidth}
\centering
\begin{tabular}{|c|c|c|}
\hline
$q_{pb}(x,x)$ & $q_{pb}(x,y)$ & $q_{pb}(y,x)$ \\ \hline
\multicolumn{3}{|c|}{$\text{Case I: For } x,y\in A, x \neq y $} \\ \hline
$x$ & $\max \{ 2x, y \}+x$ & $\max \{ 2y, x \}+y$ \\ \hline
\multicolumn{3}{|c|}{$\text{Case II: For }x,y \in B, x \neq y$} \\ \hline
$0$ & $|x-y|^{2}$ & $|y-x|^{2}$ \\ \hline
\multicolumn{3}{|c|}{$\text{Case III: For }x \in A,y \in B$} \\ \hline
$x$ & $\max \{ 2x, y \}+x$ & $3x$ \\ \hline
\multicolumn{3}{|c|}{$\text{Case IV: For }x \in B,y \in A$} \\ \hline
$x$ & $3x$ & $\max \{ 2y, x \}+y$ \\ \hline
\end{tabular}\vspace{0.3cm}
\captionof{table}{}
\end{minipage}


\noindent
\begin{minipage}{\textwidth}
\centering
\begin{tabular}{|c|c|c|c|}
\hline
$q(x,y)$ & $q(x,z)$ & $q(z,y)$ & $q(z,z)$ \\
\hline
\multicolumn{4}{|c|}{Case I: If $x,y \in A.$}\\
\hline
 \multicolumn{4}{|c|}{Subcase I: If $z \in A.$}  \\
\hline
$\max \{ 2x, y \}+x$ & $\max \{ 2x, z \}+x$ & $\max \{ 2z, y \}+z$ & $z$ \\
\hline
  \multicolumn{4}{|c|}{Subcase II: If $z \in B.$} \\
\hline
$\max \{ 2x, y \}+x$ & $\max \{ 2x, z \}+x$ & $3z$ & $z$  \\
\hline
 \multicolumn{4}{|c|}{Case II: If $x,y \in B.$}  \\
\hline
 \multicolumn{4}{|c|}{Subcase I: If $z \in A.$}\\
\hline
 $|x-y|^{2}$ & $3x$  & $\max \{ 2z, y \}+z$ & $z$  \\
\hline
\multicolumn{4}{|c|}{Subcase II: If $z \in B.$}  \\
\hline
 $|x-y|^{2}$ & $|x-z|^{2}$ & $|z-y|^{2}$ & $0$ \\
\hline
\multicolumn{4}{|c|}{$\text{Case III: If }x \in A,y \in B$.}  \\
\hline
 \multicolumn{4}{|c|}{$\text{Subcase I: If }z \in A.$}  \\
\hline
$\max \{ 2x, y \}+x$ & $\max \{ 2x, z \}+x$  & $\max \{ 2z, y \}+z$ & $z$ \\
\hline
\multicolumn{4}{|c|}{Subcase II: If $z \in B.$} \\
\hline
$\max \{ 2x, y \}+x$ &  $\max \{ 2x, z \}+x$ & $|z-y|^{2}$ & $0$\\
\hline
\multicolumn{4}{|c|}{$\text{Case IV: If }x \in B,y \in A$}  \\
\hline
  \multicolumn{4}{|c|}{$\text{Subcase I: If }z \in A.$} \\
\hline
 $3x$& $3x$ & $\max \{ 2z, y \}+z$ & $z$ \\
\hline
 \multicolumn{4}{|c|}{Subcase II: If $z \in B.$}\\
\hline $3x$ & $|x-z|^{2}$ & $|z-y|^{2}$ & $0$\\
\hline

\end{tabular}
\captionof{table}{}
\end{minipage}
Now, we define the self mappings $U$ and $V$ on $X$ as: \begin{center}
     $Ux = \begin{cases}
\dfrac{1}{6}\sin \dfrac{x}{2}, & \text{if } x \in A \\
x^{2}, & \text{ otherwise. } 
\end{cases}\quad \text{ and } \hspace{0.5cm} Vx = \begin{cases}
\dfrac{\log (x+1)}{6}, & \text{if }  x,y \in A \\
e^{x+y}, & \text{ otherwise. }
\end{cases}$
\end{center}
Also, consider the mappings $\delta,\phi: X \times X \to [0,\infty)$ defined as:\begin{center}
     $\delta(x,y) = \begin{cases}
\cos\Big(\dfrac{x+y}{4} \Big), & \text{if } x,y \in A \\
\log(x+y), & \text{ otherwise. } 
\end{cases}\quad \text{ and } \hspace{0.5cm}\phi(x,y) = \begin{cases}
\sin\Big(\dfrac{x+y}{4}\Big), & \text{if }  x,y \in A \\
e^{x+y}, & \text{ otherwise. }
\end{cases}$
\end{center}
Note that $\lim \limits_{n,m \to \infty }q_{pb}\left(5-\dfrac{1}{n},5-\dfrac{1}{m}\right)=\lim \limits_{n,m \to \infty }\Big|5-\dfrac{1}{n}-5+\dfrac{1}{m}\Big|^{2}=\Big|\dfrac{1}{n}-\dfrac{1}{m}\Big|^{2}=0.$ 
Hence, $\left(5-\dfrac{1}{n}\right)$ is a Cauchy sequence  but not convergent as $\lim \limits_{n \to \infty} \Big(5+\dfrac{1}{n}\Big)= 5 \notin X$. So, $(X,q_{pb},s)$ is not a complete  quasi-partial $b$-metric space. On the other hand, if $(x_{n})$ is a 
 $(\delta,\phi)$-Cauchy sequence in $X$,  then $(x_{n})$ is convergent in $X$. Therefore, $(X,q_{pb},s)$ is a $(\delta,\phi)$-complete quasi-partial $b$-metric space. Next, observe that for $x_{0}= \dfrac{1}{2},\epsilon =\dfrac{9}{2},$ we have $\overline{B^{l}}_{q_{pb}}\Big(\dfrac{1}{2},\dfrac{9}{2} \Big)=A.$ We will prove that $U,V$ are  locally triangular  $(\delta, \phi)$-dominated mappings on $\overline{B^{l}}_{q_{pb}}\Big(\dfrac{1}{2},\dfrac{9}{2} \Big)$.
For this, let $x \in \overline{B^{l}}_{q_{pb}}\Big(\dfrac{1}{2},\dfrac{9}{2} \Big)$. Then,  we have
$\delta(x,Ux)=\delta \Big(x,\dfrac{1}{6}\sin \dfrac{x}{2} \Big)=\cos\Big(\dfrac{x}{4}+ \dfrac{1}{24}\sin \dfrac{x}{2} \Big)\geq \phi(x,Ux)=\phi \Big(x,\dfrac{1}{6}\sin \dfrac{x}{2} \Big)=\sin \Biggr(\dfrac{x}{4}+ \dfrac{1}{24}\sin \dfrac{x}{2} \Biggr) $
and 
$\delta(x,Vx)=\delta \Big(x,\dfrac{\log (x+1)}{6} \Big)=\cos\Big(\dfrac{6x+ \log (x+1)}{24} \Big)\geq \phi(x,Vx)=\phi\Big(x,\dfrac{ \log (x+1)}{6} \Big)=\sin \Big(\dfrac{6x+ \log (x+1)}{24} \Big)$.
Hence, $U$ and $V$   are  locally   $(\delta, \phi)$-dominated mappings on $\overline{B^{l}}_{q_{pb}}\Big(\dfrac{1}{2},\dfrac{9}{2}\Big)$. Also, whenever $\delta(x,y)\geq \phi(x,y)$ and $\delta(y,z)\geq \phi(y,z)$, we have $\delta(x,z)\geq \phi(x,z).$ So, the pair $(\delta,\phi)$ is locally triangular on $\overline{B^{l}}_{q_{pb}}\Big(\dfrac{1}{2},\dfrac{9}{2} \Big)$.  Now take $\psi(x)=\dfrac{x}{6}$, for each $x \in [0,\infty)$. Clearly, $\psi \in \Psi_{s}$ as $ \sum_{i=0}^{\infty}s^{i} \psi^{i}(x)=\sum_{i=0}^{\infty}\dfrac{x}{3^{i}}$ converges for each $x \in X$. For all $j \in \mathbb{N}\cup \{0\}$, consider
\begin{align*}
  \sum_{i=0}^{j}s^{i+1} \psi^{i} \Big( \max  \Big\{ &q_{pb}\Big(\dfrac{1}{2}, U\dfrac{1}{2}\Big),q_{pb}\Big(U\dfrac{1}{2},\dfrac{1}{2}\Big)\Big\} \Big) =\\   &\sum_{i=0}^{j}2^{i+1} \psi^{i}\Big( \max  \Big\{ \max \Big\{ 1, \frac{1}{6}\sin \frac{1}{4} \Big\}+ \dfrac{1}{2}\Big\}, \max \Big\{  \frac{1}{3}\sin \frac{1}{4}  ,\dfrac{1}{2}\Big\}+ \frac{1}{6} \sin \dfrac{1}{4}\Big)  \\ &\hspace{3.9cm}= \sum_{i=0}^{j}2^{i+1} \psi^{i}\Big(\frac{3}{2}\Big)\\
  &\hspace{3.9cm}=2\Big( 1+\frac{1}{3}+\Big(\frac{1}{3}\Big)^{2}+...+\Big(\frac{1}{3}\Big)^{j}\Big)\Big(\frac{3}{2}\Big)\\ &\hspace{3.9cm}=2\sum_{i=0}^{\infty}\Big(\frac{1}{3}\Big)^{i}\Big(\frac{3}{2}\Big)\\&\hspace{3.9cm}=\epsilon.
\end{align*}
Let $x,y \in \overline{B^{l}}_{q_{pb}}\Big(\dfrac{1}{2},\dfrac{9}{2}\Big ).$  Then, $\delta(x,y) \geq \phi(x,y)$. Consider 
\begin{align}
          &\max \{ q_{pb}(Ux,Vy),q_{pb}(Vy,Ux) \}= \max \Big\{ q_{pb}\Big(\frac{1}{6}\sin \dfrac{x}{2},\frac{\log (y+1)}{6}\Big), q_{pb}\Big(\dfrac{\log (y+1)}{6},\frac{1}{6}\sin \dfrac{x}{2}\Big) \Big\}\notag\\& =\max \Big\{ \max \Big \{\frac{1}{3}\sin \dfrac{x}{2}, \frac{\log (y+1)}{6} \Big \}+ \frac{1}{6}\sin \dfrac{x}{2}, \max \Big \{ \frac{\log (y+1)}{3},\frac{1}{6}\sin \dfrac{x}{2} \Big \}+ \frac{\log (y+1)}{6} \Big \}
\end{align}
Now, we have the following two cases:\\
Case I: If $x\leq y$, then from inequality (18), we get \begin{align*}\max  \{ q_{pb}(Ux,Vy),q_{pb}(Vy,Ux) \}&\leq \frac{ \max \Biggr \{ \max \{ 2x,y\} + x, \max \{2y, \dfrac{\log (y+1)}{6} \}+y \Biggr \}}{6}\\&= \dfrac{q_{pb(y,Vy)}}{6}
\leq \psi(M_{s}(x,y)).\end{align*}
Case II: If $x> y$, then from inequality (18), we get
\begin{align*}
\hspace{-1.8cm}\max  \{ q_{pb}(Ux,Vy),q_{pb}(Vy,Ux) \} &=  \max \Biggr \{  \frac{1}{3}\sin \dfrac{x}{2}, \frac{\log (y+1)}{6}\Biggr\} + \frac{1}{6}\sin \dfrac{x}{2} \\&\leq \frac{\max \{ 2x,y\} + x}{6}\\ &\leq \dfrac{q_{pb}(x,y)}{6} \leq \psi(M_{s}(x,y)).\end{align*}
Hence, $\max  \{ q_{pb}(Ux,Vy),q_{pb}(Vy,Ux) \}\leq  \psi(M_{s}(x,y))$,  for all $x \in \overline{B^{l}}_{q_{pb}}\Big(\dfrac{1}{2},\dfrac{9}{2} \Big).$\\ Also, consider
\begin{align*}q_{pb}(Vy,Ux)&= q_{pb}\Big(\dfrac{\log (y+1)}{6},\dfrac{1}{6}\sin \dfrac{x}{2}\Big)\\&= \max \Big \{ \dfrac{\log (y+1)}{3},\dfrac{1}{6}\sin \dfrac{x}{2} \Big \}+ \dfrac{\log (y+1)}{6}\\&\leq \max \Big \{ \dfrac{y}{6} ,\dfrac{x}{6} \Big \} + \dfrac{y}{6}.\end{align*}
Now, consider the following cases:\\
 Case I: If $x\leq y$, then  we have
$q_{pb}(Vy,Ux) \leq \dfrac{y}{3} \leq \max \{ 2x, y \} +x = q_{pb}(x,y).$\\
Case II: If $x > y$, then
 we get
$q_{pb}(Vy,Ux) \leq \dfrac{x}{3} \leq \max \{ 2x, y \} +x = q_{pb}(x,y).$\\
Thus, $q_{pb}(Vy,Ux) \leq q_{pb}(x,y)$, for all $x \in \overline{B^{l}}_{q_{pb}}\Big(\dfrac{1}{2},\dfrac{9}{2} \Big).$ 
Also, $\delta(x,y)\geq \phi(x,y),$ for all $x \in \overline{B^{l}}_{q_{pb}}\Big(\dfrac{1}{2},\dfrac{9}{2} \Big).$ Thus, all the conditions of Theorem 1 are satisfied. Hence, the pair $(U,V)$ has a unique common fixed point $0$.
\end{example}
\vspace{-0.5cm}We will extend the Theorem 1 to single mapping, $\delta$-dominated mapping  and $(\delta,\phi)$-complete metric space as follows: 
\vspace{-0.4cm}\begin{corollary}
    \label{thm1}
Let $x_{0} \in X, \epsilon > 0$ and $ \delta, \phi : X \times X \to [0,\infty)$ such that the pair $(\delta, \phi)$ is locally triangular  on $\overline{B^{l}}_{q_{pb}}(x_{o}, \epsilon)$. Also, $(X,q_{pb}, s)$ be an LR-$(\delta,\phi)$ 
 complete quasi-partial $b$-metric space with coefficient $s\geq 1$ and $U: X \to X $ be  locally $(\delta, \phi)$-dominated  self map  on $\overline{B^{l}}_{q_{pb}}(x_{o}, \epsilon)$. Suppose that 
 there exists a function $\psi \in \Psi_{s}$ such that $U$  has the following properties: \\
     (i) For any $x,y \in \overline{B^{l}}_{q_{pb}}(x_{o}, \epsilon)$ such that either $\delta(x,y) \geq \phi(x,y) $ or $\delta(y,x) \geq \phi(y,x) $
     \begin{equation*}
         \text{implies } \max \{ q_{pb}(Ux,Uy),  q_{pb}(Uy,Ux) \} \leq \psi (M_{s}(x,y)),
     \end{equation*}
     where
     \begin{equation*}
        \hspace{-0.4cm} M_{s}(x,y)= \max
\Biggr \{ q_{pb}(x,y),    q_{pb}(x,Ux), q_{pb}(y,Uy), \frac{q_{pb}(x,Uy)+q_{pb}(y,Ux)-q_{pb}(x,x)}{2s}\Biggr \}.\end{equation*}
\begin{equation*}
      \text{(ii) } q_{pb}(Uy,Ux) \leq q_{pb}(x,y),\text{ for all } x,y \in \overline{B^{l}}_{q_{pb}}(x_{o}, \epsilon). \hspace{5.2cm}\end{equation*}
\begin{equation*}
    \text{(iii)}\sum_{i=0}^{j}s^{i+1} \psi^{i} ( \max  \{ q_{pb}(x_{0}, Ux_{0}), q_{pb}(Ux_{0},x_{0})\} ) \leq \epsilon, \text{ for all 
 } j \in \{ 0\}\cup \mathbb{N}.\hspace{3cm}
    \end{equation*}
  Then $U$  has atleast one  fixed point. Moreover, the fixed point  set is singleton if any pair $(x,y)$ from the set of   fixed points of   $U$ satisfies $\delta(x,y)\geq \phi(x,y).$
\end{corollary}
\vspace{-1cm}\begin{corollary}
    \label{thm1}
Let $x_{0} \in X, \epsilon > 0$ and $ \delta, \phi : X \times X \to [0,\infty)$ such that the pair $(\delta, \phi)$ is locally triangular  on $\overline{B}_{d}(x_{o}, \epsilon)$. Also, $(X,d)$ be a $(\delta,\phi)$-complete metric space and $U: X \to X $ be locally $(\delta, \phi)$-dominated  self map  on $\overline{B}_{d}(x_{o}, \epsilon)$. Suppose that 
 there exists a function $\psi \in \Psi_{s}$ such that  $U$ has the following properties:\\
     (i) For any $x,y \in \overline{B}_{d}(x_{o}, \epsilon)$ such that  
     $
      \delta(x,y) \geq \phi(x,y)     \text{ implies } d(Ux,Uy) \leq \psi (M(x,y)),
    $
     \begin{equation*}
        \hspace{-0.9cm} \text{where }M(x,y)= \max
\Biggr \{ d(x,y),    d(x,Ux),d(y,Uy), \frac{d(x,Uy)+d(y,Ux)-d(x,x)}{2}\Biggr \}.\end{equation*}
\begin{equation*}
   \text{(ii) }d(Ux,Uy)\leq d(x,y) \text{ for all } \overline{B}_{d}(x_{o}, \epsilon).\hspace{7.1cm}
\end{equation*}
\begin{equation*}
    \text{(iii) }\sum_{i=0}^{j} \psi^{i} ( d(x_{0}, Ux_{0}) ) \leq \epsilon, \text{ for all 
 } j \in \{ 0\}\cup \mathbb{N}.\hspace {6.2cm}
    \end{equation*}
Then $U$ has atleast one  fixed point. Moreover, the  fixed point  set is singleton if any pair $(x,y)$ from the set of   fixed points of   $U$ satisfies $\delta(x,y)\geq \phi(x,y).$
\end{corollary}
\vspace{-1cm}\begin{corollary}
\label{thm1}
Let $x_{0} \in X, \epsilon > 0$ and $\delta: X \times X \to [0,\infty)$ such that  $\delta$ is locally triangular  on $\overline{B^{l}}_{q_{pb}}(x_{o}, \epsilon)$. Also, $(X,q_{pb}, s)$ be an LR-$(\delta,\phi)$ 
 complete quasi-partial $b$-metric space with coefficient $s\geq 1$ and $U,V: X \to X $ be two locally $\delta$-dominated  self mappings  on $\overline{B^{l}}_{q_{pb}}(x_{o}, \epsilon)$. Suppose that 
 there exists a function $\psi \in \Psi_{s}$ such that the pair $(U,V)$ have the following properties: \\
 (i) For any $x,y \in \overline{B^{l}}_{q_{pb}}(x_{o}, \epsilon)$ such that either $\delta(x,y) \geq 1 $ or $\delta(y,x) \geq 1 $
     \begin{equation*}
         \text{implies } \max \{ q_{pb}(Ux,Vy),  q_{pb}(Vy,Ux) \} \leq \psi (M_{s}(x,y)),
     \end{equation*}
     where
     \begin{equation*}
        \hspace{-0.5cm} M_{s}(x,y)= \max
\Biggr \{ q_{pb}(x,y),    q_{pb}(x,Ux), q_{pb}(y,Vy), \frac{q_{pb}(x,Vy)+q_{pb}(y,Ux)-q_{pb}(x,x)}{2s}\Biggr \}.\end{equation*}
\begin{equation*}
       \text{(ii) }q_{pb}(Vy,Ux) \leq q_{pb}(x,y),  \text{ for all } x,y \in \overline{B^{l}}_{q_{pb}}(x_{o}, \epsilon).\hspace{7cm}\end{equation*}
 \begin{equation*}
   \text{(iii)} \sum_{i=0}^{j}s^{i+1} \psi^{i} ( \max  \{ q_{pb}(x_{0}, Ux_{0}), q_{pb}(Ux_{0},x_{0})\} ) \leq \epsilon, \text{ for all 
 } j \in \{ 0\}\cup \mathbb{N}.\hspace{4cm}
    \end{equation*}
 Then $U$ and $V$ have atleast one common fixed point. Moreover, the common fixed point  set is singleton if any pair $(x,y)$ from the set of  common fixed points of   pair $(U,V)$ satisfies $\delta(x,y)\geq \phi(x,y).$
\end{corollary}
\vspace{-0.7cm}\section{Conclusion} Common fixed point theorem on LR-$(\delta,\phi)$ quasi-partial $b$-metric spaces has been proved. This direction can be used as further extensions to  this concept to various other generalized metric spaces available in the literature.
\section*{Declarations}
\textbf{\large{Funding:}} Not applicable\vspace{0.3cm}\\
\textbf{\large{Conflict of interest:}} Authors declare no conflicts of interest\vspace{0.3cm}\\
\textbf{\large{Data availability:}} The manuscript has no associated data\vspace{0.3cm}\\
\textbf{\large{Use of AI:}} Not applicable\vspace{0.3cm}\\
\textbf{\large{Author contributions:}} All the authors contributed equally to the manuscript

\end{document}